\documentclass{llncs}
\usepackage{verbatim}
\usepackage{graphicx}
\usepackage{amsfonts}
\usepackage{amscd}
\usepackage{amssymb}
\usepackage{amsmath}

\usepackage{alltt}
\usepackage{rotating}
\usepackage{floatpag}
 \rotfloatpagestyle{empty}
\usepackage{graphicx}
\usepackage{multind}
\usepackage{times}

\usepackage{verbatim}
\usepackage{latexsym}
\usepackage{crop}
\usepackage{txfonts}
\usepackage[hyphens]{url}
\usepackage{setspace}
\usepackage{ellipsis} %

\usepackage[mathscr]{euscript} 
\usepackage{pifont} 
\usepackage[displaymath]{lineno}
\usepackage{rotating}

\usepackage{xcolor}
\usepackage{tikz}
\usetikzlibrary{calc,through,chains,shapes,arrows,trees,matrix,positioning,decorations,backgrounds,fit}

\def\tc{\hbox{:}}

\def\true{\text{true}}
\def\false{\text{false}}

\def\area{\op{area}}

\begin{document}

\title{The Strong Dodecahedral Conjecture and\\ Fejes T\'oth's Conjecture on Sphere
Packings with Kissing Number Twelve}
\author{Thomas C. Hales\thanks{{Research supported by 
NSF grant 0804189 and the Benter Foundation.}}}
\institute{University of Pittsburgh\\
\email{hales@pitt.edu}}
\maketitle

\begin{abstract}{} 
This article sketches the proofs of two theorems
  about sphere packings in Euclidean $3$-space.  The first is
  K. Bezdek's strong dodecahedral conjecture: the surface area of
  every bounded Voronoi cell in a packing of balls of radius $1$ is at
  least that of a regular dodecahedron of inradius $1$.  The second
  theorem is L. Fejes T\'oth's conjecture on sphere packings with
  kissing number twelve, which asserts that
  in $3$-space, any packing of congruent balls such that each ball is
  touched by twelve others consists of hexagonal layers.  Both proofs
  are computer assisted.  Complete proofs of these theorems appear in
  \cite{DSP} and \cite{FTK}.
\end{abstract}

     
\def\svninfo{%
  \noindent
  Document TeXed on \today. \hfill\break
  Repository Root: http://flyspeck.googlecode.com/svn \hfill\break
  SVN $LastChangedRevision: 2845 $\hfill\break
  PGF version: $\pgftypesetversion$
  }

\font\twrm=cmr8



\def\tikzfig#1#2#3{%
\begin{figure}[htb]%
  \centering
\begin{tikzpicture}#3
\end{tikzpicture}
  \caption{#2}
  \label{fig:#1}%
\end{figure}%
}

\def\tikzwrap#1#2#3#4{%
\begin{wrapfigure}{r}{#4\textwidth}
  \begin{center}
\begin{tikzpicture}#3
\end{tikzpicture}
\end{center}
  \caption{#2}
  \label{fig:#1}%
\end{wrapfigure}%
}



\renewcommand{\labelitemii}{$\circ$}
\renewcommand{\labelenumii}{\alph{enumii}}
\newenvironment{summary}
  {\begingroup\bigskip\narrower\noindent{\bf Summary.~}\it}
  {~\par\phantom{!}\endgroup\bigskip}
\def\wasitemize{\relax}
\def\uncase#1{{\sc #1}}
\def\case#1{{\sc (#1)}}
\def\firstcase#1{{\indy{Index}{named property!#1}}\relax{\sc (#1)}}
\def\claim#1{{\it  #1}}
\def\calcentry#1#2#3#4{{\smallskip{\bf #1}\quad{\tt [#2]}\quad{(#3)}\quad {#4}}} 
\def\id#1{\ensuremath{\text{\tt #1}}}

\def\indy#1#2{\index{index/#1}{#2}\relax}
\def\newterm#1{\indy{Index}{#1}{\it #1}\relax}
\def\indexterm#1{\indy{Index}{#1}{#1}\relax}
\def\fullterm#1#2{\indy{Index}{#2}{\it #1}\relax}
\def\cc#1#2{%
  \indy{Index}{computer calculation!{#1}}{\it computer calculation}
  {\footnote{\guidfoot{#1}  #2}}~\cite{website:FlyspeckProject}} 


%
\def\hide#1{}
\def\swallowed{\relax}
\def\swallow#1\swallowed{}
\newenvironment{iproved}{}{}
\newenvironment{proved}{\resetproved\begin{iproved}}{\end{iproved}}
\def\hideproof{\renewenvironment{iproved}{%
   \centerline{\it -- Proof Proofed --}
  \renewenvironment{itemize}{}{}
  \renewenvironment{enumerate}{}{}
  \def\item{\relax}
  \catcode13=12
  \swallow
}{}}
\def\showproof{\renewenvironment{iproved}{\begin{proof}}{\end{proof}}}
\def\resetproved{\if\displayallproof t\showproof\else\hideproof\fi}

\def\rating#1{\if\displayrating t%
  {{\textsc {[rating={\ensuremath {#1}}].\ }}}\else{}\fi}
\def\rz#1{\rating{#1}}
\def\cutrate{}
\def\oldrating#1{\if\displayrating t%
  {{\textsc {[former rating={\ensuremath {#1}}].\ }}}\else{}\fi}

\def\formalauthor#1{\if\verbose t{{\tt [formal proof by #1].\ }}\else{}\fi}
\def\dcg#1#2{{\if\verbose t%
  {{\tt{[DCG-#1]}}\indy{References}{ZC{#2 #1}@{DCG-#1}|page{#2}}}\else{}\fi}}
\def\tlabel#1{\label{#1}\if\verbose t{{\tt [#1].\ }%
   \indy{References}{#1|itt}}\else{}\fi}
\def\ifcverbose#1#2{\if\verbose t{{#1}}\else{#2}\fi}
\def\ifverbose#1{\ifcverbose{#1}{}}  
\def\formal#1{\relax }
\def\formaldef#1#2{\ifverbose{\texttt{[{#1} $\leftrightsquigarrow$ {#2}]}}}
\def\guid#1{\ifverbose{{\tt [#1]}}}
\def\guidfoot#1{{{\tt [#1]}}}

\setlength{\marginparwidth}{1.2in}
\def\mar#1{}

\def\dfrac#1#2{\frac{\displaystyle #1}{\displaystyle #2}}
\def\textand{\text{ \ and \ }}  

\def\emptyset{\varnothing}
\def\ups{\upsilonup} 

\def\|{\hbox{\ensuremath{\hspace{0.1em}|\hspace{-0.1em}|\hspace{0.1em}}}}
\def\mid{\ :\ }
\def\tc{\hbox{:}}
\def\norm#1#2{\hbox{\ensuremath{\|#1\unskip-\unskip{#2}\|}}}
\def\normo#1{{\|#1\|}}
\def\sland{\ \land\ }
\def\abs#1{|#1|}
\def\leftopen{(}
\def\leftclosed{[}
\def\rightopen{)}
\def\rightclosed{]}
\def\lp#1{{\llbracket{#1}\rrbracket}} 
\def\comp#1{[#1]}
\def\tangle#1{\langle #1\rangle}
\def\ceil#1{\lceil #1\rceil}
\def\floor#1{\lfloor #1\rfloor}
\def\ra{\rightsquigarrow}
\def\Ra {\Rightarrow}
\def\=#1{\accent"16 #1}
\def\ast{\ensuremath{{}^*}}

\def\CalV{{\mathcal V}}
\def\CalL{{\mathcal L}}
\def\BB{{\mathcal B}}
\def\MM{{\mathcal M}}
\def\powerset{{\mathscr P}}

\newcommand{\ring}[1]{\mathbb{#1}}
\def\F{{\mathbf F}} 

\def\v{{\mathbf v}}
\def\u{{\mathbf u}}
\def\w{{\mathbf w}}
\def\e{{\mathbf e} }  
\def\p{{\mathbf p}}
\def\q{{\mathbf q}}

\def\op#1{{\operatorname{#1}}}
\def\optt#1{{\operatorname{{\texttt{#1}}}}}

\def\atn{\op{arctan\ensuremath{_2}}}
\def\azim{\op{azim}}
\def\nd{\op{node}}
\def\sol{\operatorname{sol}}
\def\vol{\op{vol}}
\def\dih{\operatorname{dih}}
\def\Adih{\operatorname{Adih}}
\def\arc{\operatorname{arc}}
\def\rad{\operatorname{rad}}
\def\bool{\operatorname{bool}}
\def\true{\op{true}}
\def\false{\op{false}}

\def\orz{{\mathbf 0}} 
\def\Wdarto{W^0_{\text{dart}}}
\def\Wdart{W_{\text{dart}}}
\def\cell{\operatorname{cell}}
\def\dimaff{\operatorname{dim\,aff}}
\def\aff{\operatorname{aff}}
\def\card{\op{card}}

\def\mm{\hskip0.1em}
\def\nsqrt{\!\sqrt}
\def\del{\partial}
\def\doct{\delta_{oct}}
\def\dtet{\delta_{tet}}
\def\hm{{h_0}} 
\def\stab{c_{{\scriptstyle \text{stab}}}} 
\def\pqr#1{#1} 
\def\trunc#1#2{d_{#2}{#1}}

\def\ee{\varepsilonup}
\def\ocirc{}
\def\wild{{*}}  

\def\bu{{\underline{\u}}}
\def\bv{{\underline{\v}}}
\def\bw{{\underline{\w}}}
\def\bV{{\underline{V}}}

\def\smain{S_{\scriptstyle\text{main}}}

\def\svninfo{{\tt
  filename: sdodec.tex\hfill\break
  PDF generated from LaTeX sources on \today; \hfill\break
  Repository Root: https://flyspeck.googlecode.com/svn \hfill\break
  SVN $LastChangedRevision: 2986 $
  }
  }

\section{The strong hexagonal conjecture}

To describe methods, we begin with a proof of the following elementary
computer-assisted theorem in $\ring{R}^2$.

\begin{theorem} The perimeter of any bounded Voronoi cell of a packing
  of congruent balls of radius $1$ in $\ring{R}^2$ is at least
  $4\sqrt{3}$, the perimeter of a regular hexagon with inradius $1$.
\end{theorem}

If we adopt the convention that the perimeter of an unbounded Voronoi
cell is infinite, then the boundedness hypothesis can be dropped from
the theorem.

\begin{proof} Fix a bounded Voronoi cell in a packing of congruent
  balls of radius $1$ and fix a coordinate system with the center of
  the Voronoi cell at the origin.  The intersection of the Voronoi
  cell with a disk of radius $\sqrt2$ at the origin is a convex disk
  whose boundary $C$ consists of circular arcs and straight line
  segments.  The length of $C$ is no greater than the original
  perimeter of the Voronoi cell.  It suffices to show that the length
  of $C$ is at least $4\sqrt{3}$.

The boundary $C$ is determined by the set of centers $V$ of balls at
distance less than $\sqrt8$ from the origin, excluding the ball
centered at the origin.

The following piecewise linear function arises in the proof of the strong dodecahedral
conjecture in three-dimensions.  We make repeated use of it.
\begin{equation}
L(h) = \begin{cases} 
   (h_0-h)/(h_0-1),& \text{if } h \le h_0,\\
  0,&\text{otherwise},
 \end{cases}
\end{equation}
where $h_0 = 1.26$. 

Let $\u_1,\u_2\in V$ be distinct points such that
$T=\{\orz,\u_1,\u_2\}$ has circumradius less than $\sqrt2$.  Let
$\ell(\u_1,\u_2)$ be the length of the part of $C$ contained in the
convex hull of $T$, and let $\theta$ be the angle at $\orz$ between
$\u_1$ and $\u_2$ (see Figure~\ref{fig:abc}-a).The following inequality has been verified by
computer using interval arithmetic.



\begin{equation}\label{eqn:ineq} 
\ell(\u_1,\u_2) - b \theta(\u_1,\u_2) - 
c L(\normo{\u_1}/2) - c L(\normo{\u_2}/2)) \ge 0,
\end{equation}
where $b=4/3$ and $c=\sqrt3/3 - 2\pi/9 \approx -0.12$.
Equality holds when $T$ is an equilateral triangle with side $2$.

\begin{remark}
Let $\u_1(t)$ and $\u_2(t)$ be points such  that
\[
\normo{\u_1(t)}=\normo{\u_2(t)}=2,\quad \norm{\u_1(t)}{\u_2(t)}=t,
\]
that is a triangle at the origin with sides $2,2,t$.  Along this curve
with parameter $t$, the constants $b$ and $c$ are the unique choice of
constants that give the left-hand side of (\ref{eqn:ineq}) a local
minimum with value $0$ at $t=2$.
\end{remark}

\tikzfig{abc}
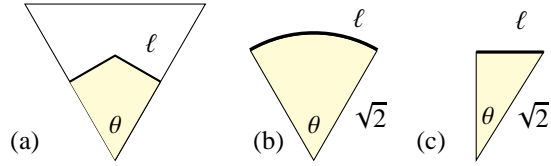
{The truncated boundary of the Voronoi cell (shaded) is partitioned into three types of pieces, indicated with a thick curve marked $\ell$.}
{
[scale=1.2]
\draw (0,0) --(1,1.732) -- (-1,1.732)--cycle;
\draw[very thick] (0.5,0.866) --(0,1.154) --(-0.5,0.866);
\draw[fill=yellow!20] (0,0)--(0.5,0.866) --(0,1.154) --(-0.5,0.866);
\path (0,0.4) node {$\theta$};
\path (0.4,1.3) node {$\ell$};
\path (-1,0.2) node {(a)};
\begin{scope}[shift={(0:2.2)}]
\draw[fill=yellow!20] (0,0) -- (60:1.4) arc(60:120:1.4) --cycle;
\draw[very thick]  (60:1.414) arc(60:120:1.414);
\path (-0.5,0.2) node {(b)};
\path (0.5,1.58) node {$\ell$};
\path (0,0.4) node {$\theta$};
\path (0.6,0.5) node {$\sqrt2$};
\end{scope}
\begin{scope}[shift={(0:4)}]
\draw[fill=yellow!20] (0,0)--(0.748,1.2)--(0,1.2)--cycle;
\draw[very thick] (0.748,1.2)--(0,1.2);
\path (0.13,0.5) node {$\theta$};
\path (0.5,1.55) node {$\ell$};
\path (-0.5,0.2) node {(c)};
\path (0.6,0.5) node {$\sqrt2$};
\end{scope}
}

The entire boundary $C$ can be partitioned into finitely many (a)
pieces lying in convex hulls of triangles $T$, (b) arcs of circles of
radius $\sqrt2$, and (c) linear segments from $\u/2$, where $\u\in V$, to a
point on the circle of radius $\sqrt2$ centered at the origin (Figure~\ref{fig:abc}).

 We extend the inequality (\ref{eqn:ineq}) to boundary arcs of type
 (b) in the form
\begin{equation}\label{eqn:ineq2}
\ell - b \theta \ge 0
\end{equation}
where $\ell$ is the length of the circular arc and $\theta$ is the
subtended angle.  This inequality is obvious, because $\ell =
\sqrt2\,\theta$, and $\sqrt2>b = 4/3$.  We extend the inequality
(\ref{eqn:ineq}) to segments of type (c) in the form
\begin{equation}\label{eqn:ineq3}
\ell(\v) - b\theta(\v) - c L(h) \ge 0,
\end{equation}
where $h =\normo{\v}/2$, $\theta$ is the subtended angle at the
origin, and $\ell = \sqrt{2-h^2}$ is the length of the the segment.
When $h\le h_0$, inequality (\ref{eqn:ineq3}) is a consequence of
inequality (\ref{eqn:ineq}), because the segment can be reflected
through a mirror to form the two segments in (\ref{eqn:ineq}).  When
$h\ge h_0$, the term $L(h)$ is zero.  In this case, basic calculus
gives the inequality.

Let $\ell_C$ be the length of $C$.  We sum the inequalities over the
boundary pieces of $C$ of types (a), (b), and (c), using inequalities
(\ref{eqn:ineq}), (\ref{eqn:ineq2}), and (\ref{eqn:ineq3}):
\begin{equation}\label{eqn:sum}
\ell_C - b (2\pi) - 2 c \sum_{\v\in V} L(\normo{\v}/2) \ge 0.
\end{equation}
The function $L$ is zero on $\{h \mid h\ge h_0\}$.  We drop such terms
from the sum.  Lemma~\ref{lemma:L} and Inequality~(\ref{eqn:sum}) give
\[
\ell_C \ge 2 \pi b + 12 c = 4\sqrt3.
\]
This proves the theorem.\qed
\end{proof}

The following lemma is used in the proof

\begin{lemma}\label{lemma:L}
  Let $V$ be a set of points contained in a closed annulus at the
  origin of inner radius $2$ and outer radius $2h_0$.  Assume that
  the mutual separation of points in $V$ is at least $2$.  Then
\[
\sum_{\v\in V} L(\normo{\v}/2) \le 6.
\]
Equality is obtained uniquely when $V$ is the set of extreme points of
a regular hexagon of circumradius $2$.
\end{lemma}

\begin{proof}
  In case $\card(V)\le 6$, by the inequality $L(\normo{\v}/2)\le 1$,
  it is clear that the sum is at most $6$, with equality uniquely
  obtained for the regular hexagon of circumradius $2$.  An easy
  estimate shows that the angles at the origin between $\u,\v\in V$ is
  greater than $\pi/4$, so that $\card(V)\le 7$.  We may therefore
  assume without loss of generality that $\card(V)=7$.

We index the $7$ points $\v_i$ by $i\in\ring{Z}/7\ring{Z}$ in their
natural cyclic order around the annulus.  Let $\theta_i$ be the angle
subtended at the origin between $\v_i$ and $\v_{i+1}$.  Let 
\[
\alpha_i = \arccos(\normo{\v}/4) - \pi/6.
\]
We have the following inequality
\begin{equation}\label{eqn:theta-alpha}
\theta_i \ge \alpha_i + \alpha_{i+1},\quad i\in \ring{Z}/7\ring{Z},
\end{equation}
which is proved by basic calculus: it follows from the intermediate
value theorem and from an explicit analytic formula for the terms in
the inequality~\cite[Lemma~6.107]{DSP}.
Further, we have
\begin{equation}\label{eqn:alpha}
  \alpha_i - 0.16\, L(\normo{\v_i}/2) - 0.32 \ge0
\end{equation}
which is also proved by basic calculus: by a second derivative test
the left-hand side is concave as a function of $\normo{\v_i}$ so that
the inequality holds if it holds at the endpoints
$\normo{\v_i}=2,2.52$, which is easily checked.

Summing $\theta_i$ over $i$ we get
\[
2\pi  = \sum_i \theta_i \ge 
2\sum_i \alpha_i \ge 2 (0.16) \sum_i L(\normo{\v_i}/2) + 14 (0.32).
\]
Computing constants, we get
\[
6 > \sum_i L(\normo{\v_i}/2).
\]
\qed
\end{proof}

\begin{remark}
  The proof can be organized in a way that carries over directly from
  two dimensions to three.  In the first step (Lemma~\ref{lemma:L}),
  $\sum L(\cdot)$ is shown to be at most the kissing number (which is $6$ in dimension two and
$12$ in dimension three).  In the
  second step, the estimate of the boundary of the truncated Voronoi
  cell is transformed into an estimate of $\sum L(\cdot)$.  The second
  step can be broken into two smaller steps: (a) use a simplex whose
  circumradius is less than $\sqrt2$ to design an inequality with a
  local minimum at the desired solution of to the Voronoi cell
  problem; (b) extend the inequality from part (a) so that it holds on
  a full geometric partition of the boundary of the truncated Voronoi
  cell.  In a final short step, sum all the inequalities to obtain the
  desired result.
\end{remark}

\section{Marchal cells}

In this section, we give details of the partition of the boundary $C$
of the truncated Voronoi cell.  The partition is based on the
partition of Euclidean space into Marchal cells~\cite{marchal:2009}.

\tikzfig{mar}
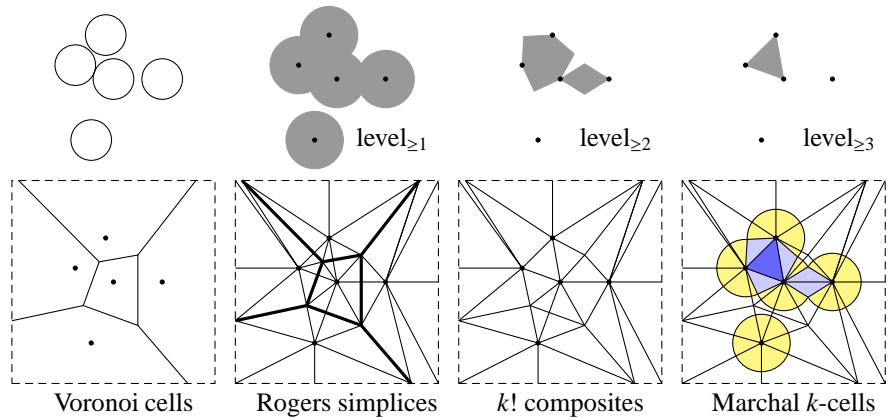
{Partitions of the plane (image source~\cite{DSP}).}
{
[scale=0.27]
\def\cf#1{\draw[fill=black] (#1) circle (0.1);}
\def\ctr{\cf{P0}\cf{P1}\cf{P2}\cf{P3}\cf{P4}}
\def\xx{5.0}
\path (0,0) node (P0) {};
\path (2.4,0) node (P1) {};
\path (100:2.2) node (P2) {};
\path (160:2) node (P3) {};
\path (250:3.2) node (P4) {};
\draw (P0) circle (1);
\draw (P1) circle (1);
\draw (P2) circle (1);
\draw (P3) circle (1);
\draw (P4) circle (1);
\begin{scope}[shift={(0:11)}]
\path (0,0) node (P0) {};
\path (2.4,0) node (P1) {};
\path (100:2.2) node (P2) {};
\path (160:2) node (P3) {};
\path (250:3.2) node (P4) {};
\def\bx#1{\draw[fill=black!40,line width=0,draw=black!40] (#1) circle (1.414);}
\bx{P0}\bx{P1}\bx{P2}\bx{P3}\bx{P4}
\path (0.5,-3)  node[anchor=west] {$\text{level}_{\ge1}$};
\ctr
\end{scope}
\begin{scope}[shift={(0:22)}]
\coordinate (P0) at (0,0);
\coordinate (P1) at (2.4,0);
\coordinate (P2) at (100:2.2);
\coordinate (P3) at (160:2);
\coordinate (P4) at (250:3.2);
\coordinate (P10) at (1.2, -0.748331);
\coordinate (P01) at  (1.2, 0.748331);
\coordinate (P02) at  (0.684303, 1.23763);
\coordinate (P20) at (-1.06633, 0.928947);
\coordinate (P30) at (-1.28171, -0.597672);
\coordinate (P03) at (-0.597672, 1.28171);
\coordinate (P23) at (-1.79446, 2.0957);
\coordinate (P32) at (-0.466949, 0.754916);
\def\bx#1#2#3#4{\draw[fill=black!40,line width=0,draw=black!40] (#1) -- (#2) -- (#3) -- (#4) -- cycle;}
\bx{P0}{P01}{P1}{P10}
\bx{P0}{P02}{P2}{P20}
\bx{P0}{P03}{P3}{P30}
\bx{P2}{P23}{P3}{P32}
\ctr
\path (0.5,-3)  node[anchor=west] {$\text{level}_{\ge2}$};
\end{scope}
\begin{scope}[shift={(0:33)}]
\coordinate (P0) at (0,0);
\coordinate (P1) at (2.4,0);
\coordinate (P2) at (100:2.2);
\coordinate (P3) at (160:2);
\coordinate (P4) at (250:3.2);
\draw[fill=black!40,draw=black!40] (P0) -- (P2) -- (P3) --(P0);
\path (0.5,-3)  node[anchor=west] {$\text{level}_{\ge3}$};
\ctr
\end{scope}
\begin{scope}[shift={(0,-10)}]
\coordinate (P0) at (0,0);
\coordinate (P1) at (2.4,0);
\coordinate (P2) at (100:2.2);
\coordinate (P3) at (160:2);
\coordinate (P4) at (250:3.2);
\coordinate (q023) at (-0.702734, 0.993058);
\coordinate (q012) at (1.2, 1.32856);
\coordinate (q014) at (1.2, -2.13945);
\coordinate (q034) at (-1.48692, -1.16149);
\coordinate (q12x) at (4.05923, 5.0);
\coordinate (q23x) at (-4.67001, 5.0);
\coordinate (q34x) at (-5., -1.90856);
\coordinate (q14x) at (3.66153, -5.);
\draw (q012)--(q023)--(q034)--(q014)--cycle;
\draw (q012)--(q12x);
\draw (q023)--(q23x);
\draw (q034)--(q34x);
\draw (q014)--(q14x);
\draw[densely dashed] (-\xx,-\xx)--(\xx,-\xx)--(\xx,\xx)--(-\xx,\xx)--cycle;
\ctr
\path (0.5,-5)  node[anchor=north] {Voronoi cells};
\end{scope}
\begin{scope}[shift={(11,-10)}]
\coordinate (P0) at (0,0);
\coordinate (P1) at (2.4,0);
\coordinate (P2) at (100:2.2);
\coordinate (P3) at (160:2);
\coordinate (P4) at (250:3.2);
\coordinate (q023) at (-0.702734, 0.993058);
\coordinate (q012) at (1.2, 1.32856);
\coordinate (q014) at (1.2, -2.13945);
\coordinate (q034) at (-1.48692, -1.16149);
\coordinate (q12x) at (4.05923, 5.0);
\coordinate (q23x) at (-4.67001, 5.0);
\coordinate (q34x) at (-5., -1.90856);
\coordinate (q14x) at (3.66153, -5.);
\coordinate (m23) at (-1.13071, 1.42531);
\coordinate (m01) at (1.2, 0);
\coordinate (m02) at (-0.191013, 1.08329);
\coordinate (m03) at (-0.939693, 0.34202);
\coordinate (m04) at (-0.547232, -1.50351);
\coordinate (m34) at (-1.48692, -1.16149);
\coordinate (e1) at (5,0);
\coordinate (e2) at (-0.382026, 5);
\coordinate (e3) at (-5,0.68404);
\coordinate (e4) at (-1.09446,-5);
\draw[very thick] (q012)--(q023)--(q034)--(q014)--cycle;
\draw[very thick] (q012)--(q12x);
\draw[very thick] (q023)--(q23x);
\draw[very thick] (q034)--(q34x);
\draw[very thick] (q014)--(q14x);
\draw (P0)--(P1);
\draw (P0)--(P2);
\draw (P0)--(P3);
\draw (P0)--(P4);
\draw (P0)--(q012);
\draw (P0)--(q023);
\draw (P0)--(q034);
\draw (P0)--(q014);
\draw (P1)--(q012)--(P2)--(P3)--(P4)--(q014)--cycle;
\draw (P2)--(q023)--(P3);
\draw (P1)--(q12x);
\draw (P1)--(e1) (P1)--(5,5) (P1)--(q12x) (P1)--(q14x) (P1)--(5,-5);
\draw (P2)--(q12x) (P2)--(e2) (P2)--(q23x);
\draw (P3)--(q23x) (P3)--(-5,5) (P3)--(e3) (P3)--(q34x);
\draw (P4)--(q34x) (P4)--(-5,-5) (P4)--(e4) (P4)--(q14x);
\draw[densely dashed] (-\xx,-\xx)--(\xx,-\xx)--(\xx,\xx)--(-\xx,\xx)--cycle;
\ctr
\path (0.5,-5)  node[anchor=north] {Rogers simplices};
\end{scope}
\begin{scope}[shift={(22,-10)}]
\coordinate (P0) at (0,0);
\coordinate (P1) at (2.4,0);
\coordinate (P2) at (100:2.2);
\coordinate (P3) at (160:2);
\coordinate (P4) at (250:3.2);
\coordinate (q023) at (-0.702734, 0.993058);
\coordinate (q012) at (1.2, 1.32856);
\coordinate (q014) at (1.2, -2.13945);
\coordinate (q034) at (-1.48692, -1.16149);
\coordinate (q12x) at (4.05923, 5.0);
\coordinate (q23x) at (-4.67001, 5.0);
\coordinate (q34x) at (-5., -1.90856);
\coordinate (q14x) at (3.66153, -5.);
\coordinate (m23) at (-1.13071, 1.42531);
\coordinate (m01) at (1.2, 0);
\coordinate (m02) at (-0.191013, 1.08329);
\coordinate (m03) at (-0.939693, 0.34202);
\coordinate (m04) at (-0.547232, -1.50351);
\coordinate (m34) at (-1.48692, -1.16149);
\coordinate (e1) at (5,0);
\coordinate (e2) at (-0.382026, 5);
\coordinate (e3) at (-5,0.68404);
\coordinate (e4) at (-1.09446,-5);
\draw (q034)--(q014);
\draw (q012)--(q12x);
\draw (q034)--(q34x);
\draw (q014)--(q14x);
\draw (P0)--(P1);
\draw (P0)--(P2);
\draw (P0)--(P3);
\draw (P0)--(P4);
\draw (P0)--(q012);
\draw (P0)--(q034);
\draw (P0)--(q014);
\draw (P1)--(q012)--(P2)--(P3)--(P4)--(q014)--cycle;
\draw (P1)--(q12x);
\draw (P1)--(e1) (P1)--(5,5) (P1)--(q12x) (P1)--(q14x) (P1)--(5,-5);
\draw (P2)--(q12x) (P2)--(e2) (P2)--(q23x);
\draw (P3)--(q23x) (P3)--(-5,5) (P3)--(e3) (P3)--(q34x);
\draw (P4)--(q34x) (P4)--(-5,-5) (P4)--(e4) (P4)--(q14x);
\draw[densely dashed] (-\xx,-\xx)--(\xx,-\xx)--(\xx,\xx)--(-\xx,\xx)--cycle;
\ctr
\path (0.5,-5)  node[anchor=north] {$k!$ composites};
\end{scope}
\begin{scope}[shift={(33,-10)}]
\coordinate (P0) at (0,0);
\coordinate (P1) at (2.4,0);
\coordinate (P2) at (100:2.2);
\coordinate (P3) at (160:2);
\coordinate (P4) at (250:3.2);
\coordinate (P10) at (1.2, -0.748331);
\coordinate (P01) at  (1.2, 0.748331);
\coordinate (P02) at  (0.684303, 1.23763);
\coordinate (P20) at (-1.06633, 0.928947);
\coordinate (P30) at (-1.28171, -0.597672);
\coordinate (P03) at (-0.597672, 1.28171);
\coordinate (P23) at (-1.79446, 2.0957);
\coordinate (P32) at (-0.466949, 0.754916);
\coordinate (q023) at (-0.702734, 0.993058);
\coordinate (q012) at (1.2, 1.32856);
\coordinate (q014) at (1.2, -2.13945);
\coordinate (q034) at (-1.48692, -1.16149);
\coordinate (q12x) at (4.05923, 5.0);
\coordinate (q23x) at (-4.67001, 5.0);
\coordinate (q34x) at (-5., -1.90856);
\coordinate (q14x) at (3.66153, -5.);
\coordinate (m23) at (-1.13071, 1.42531);
\coordinate (m01) at (1.2, 0);
\coordinate (m02) at (-0.191013, 1.08329);
\coordinate (m03) at (-0.939693, 0.34202);
\coordinate (m04) at (-0.547232, -1.50351);
\coordinate (m34) at (-1.48692, -1.16149);
\coordinate (e1) at (5,0);
\coordinate (e2) at (-0.382026, 5);
\coordinate (e3) at (-5,0.68404);
\coordinate (e4) at (-1.09446,-5);
\def\bx#1{\draw[fill=yellow!60,line width=0,draw=yellow!60] (#1) circle (1.414);}
\bx{P0}\bx{P1}\bx{P2}\bx{P3}\bx{P4}
\draw (P4) circle (1.414);
\draw (P1) circle (1.414);
\draw (P2) circle (1.414);
\draw (P3) circle (1.414);
\draw (P0) circle (1.414);
\def\bx#1#2#3#4{\draw[fill=blue!20,line width=0,draw=blue!20] (#1) -- (#2) -- (#3) -- (#4) -- cycle;}
\bx{P0}{P01}{P1}{P10}
\bx{P0}{P02}{P2}{P20}
\bx{P0}{P03}{P3}{P30}
\bx{P2}{P23}{P3}{P32}
\draw[fill=blue!60,draw=blue!60] (P0) -- (P2) -- (P3) --(P0);
\draw (q034)--(q014);
\draw (q012)--(q12x);
\draw (q034)--(q34x);
\draw (q014)--(q14x);
\draw (P0)--(P1);
\draw (P0)--(P2);
\draw (P0)--(P3);
\draw (P0)--(P4);
\draw (P0)--(q012);
\draw (P0)--(q034);
\draw (P0)--(q014);
\draw (P1)--(q012)--(P2)--(P3)--(P4)--(q014)--cycle;
\draw (P1)--(q12x);
\draw (P1)--(e1) (P1)--(5,5) (P1)--(q12x) (P1)--(q14x) (P1)--(5,-5);
\draw (P2)--(q12x) (P2)--(e2) (P2)--(q23x);
\draw (P3)--(q23x) (P3)--(-5,5) (P3)--(e3) (P3)--(q34x);
\draw (P4)--(q34x) (P4)--(-5,-5) (P4)--(e4) (P4)--(q14x);
\draw (P0)--(P01)--(P1)--(P10)--cycle;
\draw (P0)--(P02)--(P2)--(P23)--(P3)--(P30)--cycle;
\ctr
\path (0.5,-5)  node[anchor=north] {Marchal $k$-cells};
\draw[densely dashed] (-\xx,-\xx)--(\xx,-\xx)--(\xx,\xx)--(-\xx,\xx)--cycle;
\end{scope}
}

Figure~\ref{fig:mar} shows a packing $V$ of cardinality five.  We use
the constant $\sqrt2$ to partition the plane into levels $0,\ldots,3$.
Every point has level $\ge 0$.  For every $\v\in V$, we form a closed
disk of radius $\sqrt2$.  A point at level $\ge1$ is a point that lies
inside some such disk.  We form a closed rhombus of side $\sqrt2$ for
every pair of distinct points in $V$ whose separation is less than
$\sqrt8$.  By construction, the two points in $V$ are opposite
vertices of the rhombus.  A point of level $\ge2$ is a point that lies
some rhombus.  We form a closed triangle for every triple of distinct
points in $V$ whose circumradius is less than $\sqrt2$.  A point of
level $3$ is a point that lies inside some such triangle.  No point has
level $\ge4$. A point of level $k$ is a point of level $\ge k$ that
does not have level $\ge k+1$.

The points of a given level can be further partitioned using the
Rogers partition of the plane into
simplices~\cite{Rogers:1958:Packing}.  For each $k=0,\ldots,3$, Rogers
simplices that meet the set of level $k$ can be naturally grouped into
collections of $k!$ simplices.  If $P$ is the union of the $k!$
simplices, then the set of points of level $k$ in $P$ is called a {\it
  Marchal $k$-cell} $P_k$.

The construction can be generalized to three or more dimensions,
again using the parameter $\sqrt2$.   In $n$ dimensions, the levels extend
from $0$ to $n+1$ in an analogous manner.  Let $S\subset V$ be a set of
cardinality $k+1$ whose circumradius is less than $\sqrt2$.
The shapes used to define level sets are the convex hulls of
\[
S \cup X_S
\]
where $X_S$ is the set of points at equidistance $\sqrt2$ from every
point of $S$. The set $X_S$ is sphere of dimension $n-\card(S)$. 
When $n=3$, the shapes are balls of radius $\sqrt2$,
bicones, bipyramids, and tetrahedra (Figure~\ref{fig:table}).  Again
in higher dimensions, the Rogers simplices can be grouped into
collections of $k!$ simplices, giving Marchal $k$-cells $P_k$, at each
level $k$.

\begin{figure}[h!]
\centering
\begin{tabular}{l l l l l}
\hline
$\op{card}(S)$~~ &$X_S$  
&$\op{conv}(S\cup X_S)$\\ [0.5ex]
\hline \\
1 &sphere& ball\\
2& circle& bi-cone\\
3& pair of points~~~& bi-pyramid\\
4& $\emptyset$ & simplex\\
 [1ex]
\hline
\end{tabular}
\caption{Convex hulls used to construct level sets in three dimensions}
\label{fig:table}
\end{figure}

\begin{remark} Marchal introduced cells to show that the Kepler
  conjecture in three-dimensions can be reduced to an inequality of
  the form
\[
\sum_{\v\in V} M(\normo{\v}) \le 12,
\]
where $M$ is a certain quartic polynomial, and $V$ is a finite packing 
contained in a closed annulus of inner radius $2$ and outer
radius $\sqrt8$.  
\end{remark}

\begin{remark}
The book~\cite{DSP} strengthens Marchal's argument to reduce
the Kepler conjecture to the inequality
\begin{equation}\label{eqn:L12}\tag{$L_{12}$}
\sum_{\v \in V} L(\normo{\v}/2) \le 12,
\end{equation}
where $L$ is the function defined above, and $V$ is a packing in the
closed annulus with inner radius $2$ and outer radius $2h_0$.  (In
adapting this inequality from dimension two to dimension three, the
two-dimensional kissing number $6$ is replaced with the
three-dimensional kissing number $12$.) The book also gives a
computer-assisted proof of the inequality (\ref{eqn:L12}), to obtain a
new proof of the Kepler conjecture.  This article shows how to deduce
the strong dodecahedral conjecture and Fejes T\'oth's 
conjecture on packings with kissing number twelve from (\ref{eqn:L12}).
\end{remark}

\begin{remark}
An old conjecture by L. Fejes
T\'oth~\cite[p.~178]{Toth:1972:Lagerungen} asserts that the minimum of
\begin{equation}\label{eqn:FT}
\sum_{\v\in V} \normo{\v},
\end{equation}
is $24+ 14/\sqrt{27}\approx 26.69$
as $V$ runs over packings of cardinality $13$ contained in a closed
annulus with inner radius $2$ and outer radius $\sqrt8$.  The
inequality (\ref{eqn:L12}) gives the best known result:
\begin{equation}
\sum_{\v\in V} \normo{\v}\ge 24 + 2 h_0 = 26.52.
\end{equation}
The inequality (\ref{eqn:L12}) also gives upper bounds for the Tammes
problem when $\card(V)=13,14,15$, but these upper bounds are weaker
than those known by semi-definite programming~\cite{BV08}.
\end{remark}

\section{Strong Dodecahedral Conjecture}

This section sketches a proof of the strong dodecahedral
conjecture~\cite{Bezdek00}:

\begin{theorem}
  The surface area of every bounded Voronoi cell in a packing of balls
  of radius $1$ is at least the surface area of a regular dodecahedron
  of inradius $1$.  Equality is obtained exactly when the bounded
  Voronoi cell is itself a regular dodecahedron.
\end{theorem}

\begin{remark} Fejes T\'oth's classical dodecahedral
  conjecture~\cite{Toth:1943:MZ} is the same statement, replacing {\it
    surface area} with {\it volume}.  The classical dodecahedral
  conjecture was proved by McLaughlin~\cite{Hales:2010:Dodec}.  To
  deduce the volume statement from the surface area statement, it is
  enough to use the volume formula $B h/3$ for a tetrahedron, where
  $B$ is its base area (the face of a Voronoi cell), and $h\ge 1$ is
  its height.
\end{remark}

\begin{proof}
  We pick coordinates so that the center of a chosen bounded Voronoi
  cell is at the origin.  As in the two-dimensional case, we may
  intersect the Voronoi cell with a closed ball of radius $\sqrt2$.
  The boundary $C$ after truncation is no greater than before.  Let
  $V$ be the set of centers of the packing in the annulus with inner
  radius $2$ and outer radius $\sqrt8$.

There is a partition $C\cap P_k$ of  $C$ associated with the set of
Marchal $k$-cells $P_k$ associated with  Rogers simplices at the origin. 
Write 
\[
\area(C) = \sum_{P_k} \area(C\cap P_k),
\]
for the areas of the various contributions.  Write $\sol(P_k)$ for the
solid angle of the Marchal cell at the origin, and write
$\dih(P_k,\v)$ for the dihedral angle of a Marchal cell $P_k$ along
the edge through the line through $\{\orz,\v\}$.

As a reference cell, we let $P_{D,4}$ be a Marchal $4$-cell of the
packing giving the regular dodecahedron of inradius $1$.  There exist
constants $a_D$ and $b_D>0$ such that
\begin{equation}\label{eqn:dodec}
\area(C\cap P_k) +  3 a_D \sol(P_k) + 
3 b_D \sum_{\v\in P_k\cap V} L(\normo{\v}/2)\dih(P_k,\v) \ge0,
\end{equation}
for all $k$-cells $P_k$ and for all $V$.  The constants $a_D\approx -0.581$ and
$b_D\approx 0.0232$ are uniquely determined if we insist that equality
is attained when $P_k = P_{D,4}$.  This inequality has been proved by
computer by interval arithmetic.  

In more detail, the constants $a_D$
and $b_D$ are determined by a $1$-dimensional family
$P_4(t)$ of tetrahedra with sides $2,2,2,t,t,t$, for $t\in \ring{R}$
where the three edges of length $2$ meet at the origin.  When
$t=t_D\approx 2.1029$ (the separation of balls in the arrangement
giving the regular dodecahedron), $P_4(t_D)$ is congruent to
$P_{D,4}$.  Forcing the equality to be exact and the derivative to
vanish when $t=t_D$, we obtain two linear equations in two unknowns
that determine $a_D$ and $b_D$.

If we sum (\ref{eqn:dodec}) over all cells, the solid angles sum to
$4\pi$, dihedral angles sum to $2\pi$, and $L$ sums to at most $12$ by
Ineq. (\ref{eqn:L12}), giving
\[
\area(C) = \sum_{P_k}\area(C\cap P_k) \ge -3 a_D 4\pi - 3 b_D (2\pi) (12).
\]
By the choice of $a_D$ and $b_D$, equality is obtained for the boundary
$C_D$ of the regular dodecahedron, 
\[
\area(C_D) = -3 a_D 4\pi - 3 b_D (2\pi) (12).
\]
Hence $\area(C)\ge \area(C_D)$. This is the desired conclusion.  (The
circumradius of the regular dodecahedron is less than $\sqrt2$ so that
$C_D$ is both the truncated and untruncated boundary.)  Tracing
through the case of equality, inequality (\ref{eqn:dodec}) is an
equality exactly when the cell has measure zero or is congruent to
$P_{D,4}$.
\end{proof}

\section{Fejes T\'oth's Conjecture on Packings with Kissing Number Twelve}

L. Fejes T\'oth conjectured the following result in
1969~\cite{Fejes-Toth:69},~\cite{Fejes-Toth:89}. 

\begin{theorem}
  In $3$-space any packing of equal ball such that each ball is touched
  by twelve others consists of hexagonal layers.
\end{theorem}

The proof of this theorem is much longer than the proof of the strong
dodecahedral conjecture.  This section describes the proof strategy.
The details are found in~\cite{FTK}.

It is enough to prove that the contact pattern of every ball is the
hexagonal-close packing (HCP) or face-centered cubic (FCC) kissing
arrangement, because these can only be extended in hexagonal layers.
In fact, the HCP piece has a preferred plane of symmetry.  Once a
single HCP piece occurs, the preferred plane must be filled with
HCP pieces.  A plane forces another hexagonal layer above it and another
hexagonal layer below it, leading to a packing of hexagonal layers.
If no HCP piece occurs, the packing is the face-centered cubic
packing, which also consists of hexagonal layers.

\begin{lemma} Let $V$ be a packing in which every ball touches twelve
  others.  Then for all distinct $\u,\v\in V$, either
  $\norm{\u}{\v}=2$ or $\norm{\u}{\v}\ge 2h_0$.
\end{lemma}

\begin{proof}
Let $ \u_1,\ldots, \u_{12}$ be the twelve kissing points
  around $\u$.  Assume that $\v\ne \u_i,\u$.  By
  Inequality~(\ref{eqn:L12}),
\[
   L(h( \u, \v))  + 12 
  =  L(h( \u, \v)) + \sum_{i=1}^{12} L(h( \u, \u_i))  \le 12.
\]
This implies that $L(h( \u, \v))\le 0$, so $\norm{ \u}{ \v}\ge 2h_0$.\qed
\end{proof}

\subsection{graph classification problems}

\begin{definition}
  Let $S^2$ be the sphere of radius $2$, centered at $\orz$.  Let
  $\CalV$ be the set of packings $V\subset \ring{R}^3$ such that
\begin{enumerate}
\item $\card(V) = 12$,
\item $V\subset S^2$,
\item $\norm{\u}{\v} \in \{0,2\}\cup
  \leftclosed2.52,4\rightclosed$ for all $\u,\v\in V$.
\end{enumerate}
For each $V\in \CalV$, let $E_{ctc}$ be the contact graph on vertex set $V$;
that is, the set of $\{\u,\v\}\subset V$ such that $\norm{\u}{\v}=2$.
\end{definition}

Fejes T\'oth's conjecture follows from the
Inequality~(\ref{eqn:L12}), together with a proof that the
classification of graphs $(V,E_{ctc})$ with $V\in \CalV$ up to
isomorphism contains exactly two graphs: the FCC contact graph and the
HCP contact graph.

We formulate Inequality~(\ref{eqn:L12}) as a graph classification
problem as well.  The inequality holds trivially for a finite packing
of cardinality at most $12$.  For a contradiction, we may assume that
$V$ belongs to the set of finite packings of cardinality at least $13$,
contained in a closed annulus of radii $[2,2h_0]$ and that violate the
inequality:
\[
\sum_{\v\in V} L(\normo{v}/2) > 12.
\]
Let $E_{std}$ be the set of edges $\{\u,\v\}\subset V$ such that
$2\le\norm{\u}{\v}\le 2h_0$.  The graph classification problem
equivalent to ($L_{12}$) is that the set of graphs $(V,E_{std})$, with
$V$ from this set of counterexamples, is empty.

In summary, the proof of Fejes T\'oth's  conjecture consists of
two graph classification problems: one for the contact graphs
$(V,E_{ctc})$ involving vertex sets of cardinality $12$ and one for
the graphs $(V,E_{std})$ involving vertex sets of cardinality at least
$13$ for the Inequality~(\ref{eqn:L12}).  The proofs of these two
classification results differ in detail, but the high-level structure
is the same in both cases.  The graphs are first represented
combinatorially as hypermaps.  (A hypermap is a finite set $D$
together with three permutations $e,n,f$ of $D$ that satisfy $e n f =
I$, the trivial permutation.)  A computer program classifies the
hypermaps satisfying given combinatorial properties obtained from the
constraints imposed on the graphs $(V,E)$.  Linear programs eliminate
the extraneous cases; namely, those hypermaps that exist
combinatorially but that do not admit a geometric realization.
Finally, the inequalities used in the linear programs are proved by
computer.

\subsection{hypermap classification by computer}

The computer program that classifies hypermaps has been the subject of
a exhaustive computer code formal verification project by G. Bauer and
T. Nipkow~\cite{Nipkow:2005:Tame}.  The original scope of the project
was the set of graphs from the 1998 proof of the Kepler
conjecture~\cite{Hales:2006:DCG}, but in 2010, Nipkow extended this
work to include the classification of hypermaps needed for for the
$L_{12}$ inequality.  \footnote{There are about 25,000 graphs that
  arise in the $L_{12}$ classification and only $8$ graphs that arise
  in the contact graph.  Because of the vast difference in complexity
  of these two classification problems, our discussion will focus on
  the $L_{12}$ classification.}

\subsection{linear programs}

As mentioned above, linear programs eliminate the extraneous cases.
The technology related to the linear programming has been
significantly improved in the years following the proof of the Kepler
conjecture.  The thesis of S. Obua implemented
the formal verification of linear programming proof certificates and
used this to eliminate about $92\%$ of the graphs that appear in the
original proof of the Kepler conjecture~\cite{Obua:2008:Thesis}.  More
recent work by Solvyev has optimized the formal verification of linear
programs to such a degree that the speed of the formal verification of
a linear program rivals the speed of the unverified execution of a
linear program~\cite{Solovyev:LP}.  Work in progress by Solovyev
intends to make a full formal verification of all linear programs needed to
prove Inequality~(\ref{eqn:L12}).

The linear programs are generated in GLPK from an AMPL model that is
independent of the hypermap.  An OCaml program generates a
separate AMPL data file for each linear program.  When a single linear
program fails to eliminate a hypermap, branch and bound methods are
used to iteratively subdivide the domain into smaller pieces until
linear programs are obtained that eliminate the hypermap.  The process
that was used to obtain a system of linear programming inequalities
that works uniformly on all hypermaps was fully
automated~\cite{Hales:2010:lin-prog}.  In brief synopsis, when a
linear program fails to eliminate a hypermap, two models of
corresponding metric graph are compared, one based purely on the
linear programming estimates of lengths and angles, and a second
nonlinear model based on nonlinear relations between lengths and
angles.  A comparison of models is used to determine inadequacies in
the linearization.  This data is fed to a Mathematica program based on
various heuristics to construct a candidate nonlinear inequality.  The
inequality is then shipped to the nonlinear optimization package
{CFSQP} for extensive nonrigorous testing.  From there, a formal
specification of the inequality is generated in the proof assistant
HOL Light.  The formal specification is exported to program that uses
interval arithmetic to verify inequalities by computer; and finally,
the AMPL model is automatically updated with the new inequality.  This
process works remarkably well in practice to develop a small set of
inequalities\footnote{About 500 inequalities occur.} that
uniformly eliminate all undesired hypermaps.

\subsection{The classification of contact graphs}

As mentioned in an earlier footnote, there are about 25 thousand graphs
that arise in the $L_{12}$ classification and only $8$ graphs that
arise in the contact graph classification. For the $8$ graphs, it was
not necessary to follow the lengthy linear programming procedure
described in the previous subsection.  This final subsection sketches a
much simpler procedure to eliminate the unwanted cases.  (Two of the
eight possibilities are the HCP and FCC, and the other six cases are
unwanted.)

Five of the six are eliminated with linear programming inequalities.
The linear programs are based on the following simple inequalities:
\begin{enumerate}
\item The angles around each node sum to $2\pi$.
\item The angle of each triangle in the contact graph  equals $\arccos(1/3)\approx 1.23096$.
\item The opposite angles of each rhombus are equal.
\item Each angle of every rhombus is between $1.6292$ and $2.16672$.
\end{enumerate}

The final case is the graph shown in Figure~\ref{fig:fthex}.  It is
eliminated with the following observations.  The perimeter of a
spherical hexagon with sides $\pi/3$ is $2\pi$.  However, the hexagons in the graph are
spherically convex, and $2\pi$ is a strict upper bound on the perimeter of a
spherically convex hexagon.  Thus, this case does not admit a geometric realization
as a contact graph.   Fejes T\'oth's 
conjecture on sphere packings with kissing number twelve ensues.\qed

\tikzfig{fthex}
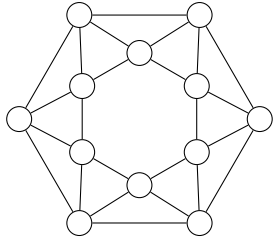
{This planar graph is not a contact graph.}
{
[scale=0.004]
\tikzstyle{every node}=[draw,shape=circle];
\path ( 400,0) node (P0) {};
\path (60:400)  node (P1) {};
\path (120:400) node (P2) {};
\path ( -400,0) node (P3) {};
\path ( -200,-346) node (P4) {};
\path ( 200,-346) node (P5) {};
\path (330:220) node (P6) {};
\path (270:220) node (P7) {};
\path(210:220) node (P8) {};
\path (30:220) node (P11) {};
\path (150:220) node (P9) {};
\path (90:220) node (P10) {}; 
\draw
  (P0) -- (P5)
  (P0) -- (P1)
  (P0) -- (P11)
  (P0) -- (P6)
  (P1) -- (P2)
  (P1) -- (P10)
  (P1) -- (P11)
  (P2) -- (P3)
  (P2) -- (P9)
  (P2) -- (P10)
  (P3) -- (P4)
  (P3) -- (P8)
  (P3) -- (P9)
  (P4) -- (P5)
  (P4) -- (P7)
  (P4) -- (P8)
  (P5) -- (P6)
  (P5) -- (P7)
  (P6) -- (P11)
  (P6) -- (P7)
  (P7) -- (P8)
  (P8) -- (P9)
  (P9) -- (P10)
  (P10) -- (P11)
;
}

\raggedright
\bibliographystyle{amsalpha} 
\bibliography{/Users/thomashales/Desktop/googlecode/flyspeck/latex/bibliography/all}

\bigskip
\noindent
\smallskip

\end{document}